\documentclass[a4paper,12pt]{article}

\usepackage{amsfonts}
\usepackage{amscd,color}
\usepackage{amsmath,amsfonts,amssymb,amscd}
\usepackage{booktabs}
\usepackage{mathrsfs}
\usepackage{indentfirst,graphicx,epsfig}
\usepackage{graphicx}
\input{epsf}
\usepackage{epstopdf}
\usepackage{caption}
\usepackage{ifpdf}
\usepackage{float}
\usepackage{footnote}
\usepackage{fullpage}
\usepackage{epsfig}
\usepackage{latexsym}
\usepackage{caption}
\usepackage{amsmath,amsthm}
\usepackage{amsfonts}
\usepackage{amssymb}
\usepackage{graphicx}
\usepackage{graphics}
\usepackage{bm}
\usepackage{color}
\usepackage{subfigure}
\usepackage{threeparttable}
\graphicspath{{\figs}}

\makeatletter

\newcommand{\Rmnum}[1]{\expandafter\@slowromancap\romannumeral #1@}
\makeatother
\makeatletter
\let\@fnsymbol\@arabic
\makeatother

\begin{document}
\newtheorem{theorem}{Theorem}[section]
\newtheorem{observation}[theorem]{Observation}
\newtheorem{corollary}[theorem]{Corollary}
\newtheorem{algorithm}[theorem]{Algorithm}
\newtheorem{problem}[theorem]{Problem}
\newtheorem{question}[theorem]{Question}
\newtheorem{lemma}[theorem]{Lemma}
\newtheorem{proposition}[theorem]{Proposition}

\newtheorem{definition}[theorem]{Definition}
\newtheorem{guess}[theorem]{Conjecture}
\newtheorem{claim}[theorem]{Claim}
\newtheorem{example}[theorem]{Example}
\makeatletter
  \newcommand\figcaption{\def\@captype{figure}\caption}
  \newcommand\tabcaption{\def\@captype{table}\caption}
\makeatother

\newtheorem{acknowledgement}[theorem]{Acknowledgement}

\newtheorem{axiom}[theorem]{Axiom}
\newtheorem{case}[theorem]{Case}
\newtheorem{conclusion}[theorem]{Conclusion}

\newtheorem{conjecture}[theorem]{Conjecture}
\newtheorem{criterion}[theorem]{Criterion}
\newtheorem{exercise}[theorem]{Exercise}
\newtheorem{notation}[theorem]{Notation}
\newtheorem{solution}[theorem]{Solution}
\newtheorem{summary}[theorem]{Summary}
\newtheorem{fact}[theorem]{Fact}
\newtheorem{remark}[theorem]{Remark}

\newcommand{\pp}{{\it p.}}
\newcommand{\de}{\em}
\newcommand{\mad}{\rm mad}

\newcommand*{\QEDA}{\hfill\ensuremath{\blacksquare}}
\newcommand*{\QEDB}{\hfill\ensuremath{\square}}

\newcommand{\qf}{Q({\cal F},s)}
\newcommand{\qff}{Q({\cal F}',s)}
\newcommand{\qfff}{Q({\cal F}'',s)}
\newcommand{\f}{{\cal F}}
\newcommand{\ff}{{\cal F}'}
\newcommand{\fff}{{\cal F}''}
\newcommand{\fs}{{\cal F},s}
\newcommand{\g}{\gamma}
\newcommand{\wrt}{with respect to }

\title{Asymptotic values of four Laplacian-type energies for matrices with degree-distance-based entries of random graphs}

\author{\small Xueliang Li$^1$, Yiyang Li$^2$, Zhiqian Wang$^1$\\
\small $^{1}$Center for Combinatorics and LPMC\\
\small Nankai University, Tianjin 300071, China\\
\small Email: lxl@nankai.edu.cn; 1522686578@qq.com\\
\small $^{2}$International Institute\\
\small China Construction Bank\\
\small Xicheng District, Beijing 100033, China\\
\small Email: liyiyang@163.com
}
\baselineskip=18pt

\date{ }
\maketitle

\begin{abstract}
Let $f(D(i, j), d_i, d_j)$ be a real function symmetric in $i$ and $j$ with the property that $f(d, (1+o(1))np, (1+o(1))np)=(1+o(1))f(d, np, np)$ for $d=1,2$. Let $G$ be a graph,
$d_i$ denote the degree of a vertex $i$ of $G$ and $D(i, j)$ denote the distance between vertices $i$ and $j$ in $G$. In this paper,
we define the $f$-weighted Laplacian matrix for random graphs in the Erd$\ddot{o}$s-R$\acute{e}$nyi random graph model $\mathcal{G}_{n, p}$, where $p\in (0, 1)$ is fixed.
Four weighted Laplacian type energies: the weighted Laplacian energy $\mathscr{LE}_f(G)$, weighted signless Laplacian energy $\mathscr{LE}^{+}_f(G)$, weighted
incidence energy $\mathscr{IE}_f(G)$ and the weighted Laplacian-energy like invariant $\mathscr{LEL}_f(G)$ are introduced and studied. We obtain the asymptotic values of $\mathscr{IE}_f(G)$
and $\mathscr{LEL}_f(G)$, and the values of $\mathscr{LE}_f(G)$ and $\mathscr{LE}_f^{+}(G)$ under the condition that $f(D(i, j), d_i, d_j)$ is a function dependent only
on $D(i, j)$. As a consequence, we get that for almost all graphs $G_p\in \mathcal{G}_{n, p}$, the energy for the matrix with degree-distance-based entries of $G_p$,
$\mathscr{E}(W_f(G_p)) < \mathscr{LE}_f(G_p),$ the Laplacian energy of the matrix, which is a generalization of a conjecture by Gutman et al.  \\[3mm]
\end{abstract}

\section[Introduction]{Introduction}

Throughout this paper, $G$ denotes a simple graph with $n$ vertices and $m$ edges. In chemistry, there is a closed relation between the molecular orbital energy levels of $\pi$-electrons in conjugated hydrocarbons and the eigenvalues of the corresponding molecular graph. For the H$\ddot{u}$chkel molecular orbital approximation, the total $\pi$-electron energy in a conjugated hydrocarbon is given by the sum of absolute values of the eigenvalues corresponding to the molecular graph $G$ in which the maximum degree is no more than 4 in general. In 1970s, Gutman in \cite{gut} extended the concept of energy to all simple graphs $G$. More on the energy
of graphs can be found in \cite{lsg}. We use $A(G)$ to represent the adjacency matrix of $G$ with spectrum $\rho_1\geq\rho_2\geq\ldots\geq\rho_n$, which are said to be the eigenvalues of graph $G$. The energy of graph $G$ is defined as
\begin{displaymath}
\mathscr{E}(G)=\sum\limits_{i=1}^{n} |\rho_i|.
\end{displaymath}

Moreover, several other energy-like quantities were introduced later, such as the Laplacian energy, signless Laplacian energy,
Laplacian-energy like invariant and the incidence energy, which have been widely studied in the mathematical and mathematical-chemical
literatures; see below for definitions.\\

The \emph{Laplacian matrix} of $G$ is defined as $L(G)=D(G)-A(G)$, where $D(G)$ is the diagonal matrix in which every diagonal entry
is equal to the degree of the corresponding vertex. Supposing $\lambda_1(L)\geq\lambda_2(L)\geq\ldots\lambda_n(L)$ are the eigenvalues
of $L(G)$, Gutman and Zhou \cite{gutzhou} defined the \emph{Laplacian energy} of $G$ as\\
\begin{displaymath}
\mathscr{LE}(G)=\sum\limits_{i=1}^{n} |\lambda_i(L)-\frac{2m}{n}|.
\end{displaymath}
The \emph{signless Laplacian energy} is defined as $L^+(G)=D(G)+A(G)$. Supposing that $\lambda_1(L^+)\geq\lambda_2(L^+)\geq\ldots\lambda_n(L^+)$ \
are the eigenvalues of $L^+$, So \emph{et al.} \cite{liu} defined the \emph{signless Laplacian energy} of $G$ as
\begin{displaymath}
\mathscr{LE}^+(G)=\sum\limits_{i=1}^{n} |\lambda_i(L^+)-\frac{2m}{n}|.
\end{displaymath}
Furthermore, Liu \emph{et al.} \cite{liu} proposed a \emph{Laplacian-energy like invariant}, which is defined as
\begin{displaymath}
\mathscr{LEL}(G)=\sum\limits_{i=1}^{n} \sqrt{\lambda_i(L)}.
\end{displaymath}
Gutman \emph{et al.} pointed out in \cite{gzf} that $\mathscr{LEL}(G)$ is more similar to $\mathscr{E}(G)$ than to $\mathscr{LE}(G)$.

Moreover, Jooyandeh \emph{et al.} \cite{joo} introduced the \emph{incidence energy} $\mathscr{IE}(G)$ of $G$,
which is defined as the sum of the singular values of the incidence matrix of $G$. Gutman \emph{et al.} showed that\\
\begin{displaymath}
\mathscr{IE}(G)=\sum\limits_{i=1}^{n} \sqrt{\lambda_i(L^+)}.
\end{displaymath}

We call these four energies Laplacian type energies since they are defined on the basis of Laplacian matrix
rather than adjacency matrix. It is not difficult to calculate these energies for a concrete graph,
just by computing the eigenvalues of $L(G)$ or $L^+(G)$. But this is impractical when $n$ is getting very large.
Some upper and lower bounds were established before, and for details we refer the reader to \cite{gkmz, rli}.
Those inequalities, however, have a common flaw that only a few specific graphs attain the equalities of the bounds.
So, they can hardly characterize the asymptotic tendency of these energies as $n\rightarrow\infty$.
In order to see clearly the asymptotic tendency, one usually employs random graphs in the so-called Erd$\ddot{o}$s-R$\acute{e}$nyi random graph model $\mathcal{G}_{n, p}$.
Recall that $\mathcal{G}_{n, p}$ consists of all graphs on $n$ vertices in which the edges are chosen independently and equally
with probability $p$, where $p\in (0, 1)$ is a constant. To study the property of general graphs, random graphs in $\mathcal{G}_{n, p}$
are very suitable objects to serve this purpose.

At first, we recall some results on random matrices, details of which can be found in \cite{wig1, wig2}.
In 1950s, Wigner studied the limiting spectral distribution of a type of random matrices, named as \emph{Wigner matrix},
denoted by $X=\{x_{ij}\}_{i, j=1}^{n}$, which satisfies the following conditions:\\
(i) \ $x_{ij} \ (i\neq j)$ are \emph{i.i.d.} random variables with variance $\sigma^{2}$, and $x_{ij}=x_{ji}$;\\
(ii) \ $x_{ii}$ are \emph{i.i.d.} random variables without any moment requirement.\\

The \emph{empirical spectral distribution (ESD)} of $X$ is defined by $\Phi_{X}(x)=\frac{1}{n}\cdot\sharp\{\lambda_{i} \ | \ \lambda_{i}\leq x, i=1, 2, \ldots, n\}$.
Then, the energy of matrix $X$ is $\mathscr{E}(X)=n\cdot\int|x|\,d\Phi_{X}(x)$.
Wigner calculated the limiting spectral distribution (LSD for short) of $X$ and obtained his famous {\it semicircle law}; see \cite{wig1, wig2}.

Define a random matrix $M_n=X_n-D_n$ to be a \emph{Markov matrix}, where $X_n$ is a Wigner matrix with all diagonal entries 0,
and $D_n$ is a diagonal matrix in which $D_n(i, i)=\sum\limits_{j\neq i} x_{ij}, i=1,\ldots n$. Bryc \emph{et al.} obtained the LSD of a Markov matrix;
see below.\\
\begin{theorem} \label{markov}
(see Bryc et al. \cite{bryc}.) Let $M_n$ be a Markov matrix such that $x_{12}$ is of expectation 0 and variance 1. Then with probability 1, $\Phi_{n^{-1/2}M_n}(x)$
converges weakly to a distribution $\Psi(x)$ as $n\rightarrow\infty$, where $\Psi$ is the free convolution of the standard semicircle law $\phi(x)=\frac{1}{2\pi}\sqrt{4-x^{2}}$
and the standard normal measure.
\end{theorem}

On the basis of the fact that the adjacency matrix $A(G_p)$ of a random graph $G_p \in \mathcal{G}_{n, p}$ is a Wigner matrix and
that $L(G_p)$ is a Markov matrix, in \cite{dll}, Du, Li and Li applied the semicircle law and Theorem \ref{markov}, and they obtained
the following results.
\begin{theorem}
For almost all random graphs $G_p \in \mathcal{G}_{n, p}$,
\begin{displaymath}
\mathscr{E}(G_n(p))=(\frac{8}{3\pi}\sqrt{p(1-p)}+o(1))\cdot n^{3/2}.
\end{displaymath}
\end{theorem}
\begin{theorem} \label{lap}
For almost all random graphs $G_p \in \mathcal{G}_{n, p}$,
\begin{displaymath}
(\frac{2\sqrt{2}}{3}\sqrt{p(1-p)}+o(1))\cdot n^{3/2}\leq \mathscr{LE}(G_n(p))\leq(\sqrt{2}\cdot\sqrt{p(1-p)}+o(1))\cdot n^{3/2}.
\end{displaymath}
\end{theorem}

In \cite{dll2}, they calculated the other three kinds of energies; see below.
\begin{theorem}\label{sign}
For almost all random graphs $G_p \in \mathcal{G}_{n, p}$,
\begin{displaymath}
[(\frac{16}{3\pi}-\sqrt{2})\cdot\sqrt{p(1-p)}+o(1)]\cdot n^{3/2}\leq\mathscr{LE}^+(G_p)\leq[(\frac{16}{3\pi}+\sqrt{2})\cdot\sqrt{p(1-p)}+o(1)]\cdot n^{3/2}.
\end{displaymath}
\end{theorem}
\begin{theorem}
For almost all random graphs $G_p \in \mathcal{G}_{n, p}$,
\begin{displaymath}
\mathscr{LEL}(G_p)=(\sqrt{p}+o(1))\cdot n^{3/2}.
\end{displaymath}
\end{theorem}
\begin{theorem}
For almost all random graphs $G_p \in \mathcal{G}_{n, p}$,
\begin{displaymath}
\mathscr{IE}(G_p)=(\sqrt{p}+o(1))\cdot n^{3/2}.
\end{displaymath}
\end{theorem}

In practical requirements of molecular chemistry, many objects are concerned with degrees of vertices and/or distances between pairs of vertices in a graph,
and many interesting matrices with entries from degrees and/or distances have been introduced; see \cite{XL, dgmmf}. They are categorized as the following three kinds. One kind of such
matrices comes from degree-based topological indices of chemical use, such as the Zagreb matrix, ABC-matrix and Harmonic matrix. Another kind of such matrices
comes from distance-based topological indices, such as the distance matrix, Harary matrix and the reverse Wiener matrix, etc. In 1994, Dobrynin and Kochetova put
forward a new topological index determined by the values of both distances between pairs of vertices and degrees of vertices, and recently this new type of indices becomes more and more
popular. We refer them as the degree-distance-based indices. From this kind of indices, it is natural to define a new kind of matrices with
mixed degree-distance-based entries, since a 2-dimensional matrix contains much more structural data than a single value of index. As one can see below,
it is essentially a distance matrix. So, we call it the \emph{weighted distance matrix} of a graph $G$. We use $d_i$ and $D(i, j)$ to represent the degree of a
vertex $i$ of $G$ and the distance between two vertices $i$ and $j$ in $G$, respectively. The definition of the matrix is given as follows.
\begin{definition}
Let $G$ be a graph, and $f(D(i, j), d_i, d_j)$ be a real function symmetric in $i$ and $j$.
The weighted distance matrix $W_f(G)$ is defined as follows: the $ij$-entry of $W_f(G)$
\begin{displaymath}
W_f(G)(i, j)=\begin{cases}
f(D(i, j), d_i, d_j), &\quad \mbox{$i\neq j$}\\
0,&\quad \mbox{$i=j$}
\end{cases}
\end{displaymath}
\end{definition}
Similarly, the \emph{weighted distance Laplacian matrix} are defined as $L_f(G)=D_f(G)-W_f(G)$, where $D_f(G)$ is a diagonal matrix in which the $ii$-entry
equals $\sum\limits_{j\neq i} f(D(i, j), d_i, d_j)$. The \emph{weighted Laplacian energy} $\mathscr{LE}_f$ and \emph{weighted signless Laplacian energy}
$\mathscr{LE}^{+}_f$ are also defined similarly, just to replace $\frac{2m}{n}$ with
$$\frac{\sum\limits_{i\neq j}f(D(i, j), d_i, d_j)}{n}.$$
If we replace $\lambda_i$ with $|\lambda_i|$, we get the definitions of \emph{weighted Laplacian-energy like invariant} $\mathscr{LEL}_f$ and \emph{weighted incidence energy} $\mathscr{IE}_f$.

Unfortunately, when we consider $W_f(G_p)$ for random graphs $G_p$, a big problem arises. This matrix is no longer a Wigner matrix since the random variables
$f(D(i, j), d_{i}, d_{j})$ (${i,j}\in \{1,2,\ldots,n\}$) are not independent. This limits the use of this method, and most previous results become unavailable.
Nevertheless, we still managed to obtain the asymptotic values of weighted energies by doing some approximation of the classical results. Throughout the paper,
we suppose that $f(1, d_i, d_j)$ and $f(2, d_i, d_j)$ satisfy the property that $f(1,(1+o(1))np, (1+o(1))np)=(1+o(1))f(1, np, np)$ and
$f(2,(1+o(1))np, (1+o(1))np)=(1+o(1))f(2, np, np)$. We give the asymptotic values of $\mathscr{LEL}_f(G_p)$ and $\mathscr{IE}_f(G_p)$, and $\mathscr{LE}_f(G_p)$ and $\mathscr{LE}_f^+(G_p)$ under the condition
that $f(D(i, j), d_i, d_j)$ is only determined by $D(i, j)$. Our main results are stated as follows.

\begin{theorem} \label{main1}
Let $f(1, d_i, d_j)$, $f(2, d_i, d_j)$ be symmetric functions satisfying that $f(1,(1+o(1))np, (1+o(1))np)=(1+o(1))f(1, np, np)$ and
$f(2,(1+o(1))np, (1+o(1))np)=(1+o(1))f(2, np, np)$. Then for almost all graphs $G_p\in \mathcal{G}_{n, p}$,\\
(i) \ if $f(1, np, np)/f(2, np, np)\rightarrow \pm\infty$, then
\begin{displaymath}
\mathscr{LEL}_f(G_p)=\sqrt{|f(1, np, np)|}(\sqrt{p}+o(1))\cdot n^{3/2},
\end{displaymath}
\begin{displaymath}
\mathscr{IE}_f(G_p)=\sqrt{|f(1, np, np)|}(\sqrt{p}+o(1))\cdot n^{3/2};
\end{displaymath}
(ii) \ if $f(1, np, np)/f(2, np, np)\rightarrow C\in\mathbf{R}$, then
\begin{displaymath}
\mathscr{LEL}_f(G_p)=\sqrt{|f(2, np, np)|}(\sqrt{|1+(C-1)p|}+o(1))\cdot n^{3/2},
\end{displaymath}
\begin{displaymath}
\mathscr{IE}_f(G_p)=\sqrt{|f(2, np, np)|}(\sqrt{|1+(C-1)p|}+o(1))\cdot n^{3/2}.
\end{displaymath}
\end{theorem}

\begin{theorem}\label{main2}
Assume that $D_1=f(1, d_i, d_j)$ and $D_2=f(2, d_i, d_j)$ are constants. Then, for almost all graphs $G_p\in \mathcal{G}_{n, p}$, we have
\begin{displaymath}
|D_1-D_2|(\frac{2\sqrt{2}}{3}\sqrt{p(1-p)}+o(1))\cdot n^{3/2}\leq \mathscr{LE}_f(G_p)
\end{displaymath}
\begin{displaymath}
\leq |D_1-D_2|(\sqrt{2}\cdot\sqrt{p(1-p)}+o(1))\cdot n^{3/2},
\end{displaymath}
and the signless Laplacian energy enjoys the inequalities
\begin{displaymath}
|D_1-D_2|[(\frac{16}{3\pi}-\sqrt{2})\cdot\sqrt{p(1-p)}+o(1)]\cdot n^{3/2}\leq
\mathscr{LE}_f^+(G_p)
\end{displaymath}
\begin{displaymath}
\leq|D_1-D_2|[(\frac{16}{3\pi}+\sqrt{2})\cdot\sqrt{p(1-p)}+o(1)]\cdot n^{3/2}.
\end{displaymath}
\end{theorem}

For the asymptotic value of adjacency energy of random graphs with degree-based weights $f(d_{i}, d_{j})$, Li, Li and Song in \cite{lls} obtained the following result.

\begin{theorem}\label{adj}
Let $f(x, y)$ be a symmetric real function. Denote by $A_f(G_p)$ the adjacency matrix of a random graph $G_p$ weighted by a degree-based function $f(d_i,d_j)$. If the function $f$ satisfies that conditions  that $|f(d_{i}, d_{j})|\leq Cn^{m}$ for some constants $C, m > 0$, and $f((1+o(1))np, (1+o(1))np)=(1+o(1))f(np, np)$ where $p\in(0, 1)$ is any fixed and independent of $n$, then for almost all graphs $G_p$ in $\mathcal{G}_{n, p}$,\\
\begin{displaymath}
\mathscr{E}(A_f(G_p)=|f(np, np)|(\frac{8}{3\pi}\sqrt{p(1-p)}+o(1))\cdot n^{3/2} \quad a.s.
\end{displaymath}
\end{theorem}

Whereas, for the asymptotic value of energy of weighted distance matrix $W_f(G_p)$ of graphs $G_p\in\mathcal{G}_{n, p}$, Li, Li and Wang in \cite{llw} recently obtained the following result.\\

\begin{theorem} \label{dis}
Let $f(D(i, j), d_{i}, d_{j})$ be a function symmetric in $i$ and $j$ satisfying the condition that $f(D(i, j), (1+o(1))np, (1+o(1))np)=(1+o(1))f(D(i, j), np, np)$. Then for almost all graphs $G_p\in \mathcal{G}_{n, p}$,
$$
\mathscr{E}(W_f(G_p))=\{(\frac{8}{3\pi}\sqrt{p(1-p)}+o(1))\cdot|D_1- D_2|+o(|f(2, np, np)|)\}\cdot n^{3/2}\quad a.s.
$$
That is, if $f(1, np, np)/f(2, np, np)\nrightarrow 1$,\\
$$
\mathscr{E}(W_f(G_p))=|D_1-D_2|(\frac{8}{3\pi}\sqrt{p(1-p)}+o(1))\cdot n^{3/2}\quad a.s.
$$
and if $f(1, np, np)/f(2, np, np)\rightarrow 1$,\\
$$
\mathscr{E}(W_f(G_p))=o(1)|f(2, np, np)|\cdot n^{3/2}\quad a.s.
$$
\end{theorem}

\begin{remark}
Note that if one considers only adjacency but not distance matrices, then by setting $f(2,d_j,d_j)=0$, it is easy to
check that $W_f(G_p)=A_f(G_p)$ and therefore Theorem \ref{dis} implies Theorem \ref{adj}. So, Theorem \ref{adj} holds if one deletes the condition
that $|f(d_{i}, d_{j})|\leq Cn^{m}$ for some constants $C, m > 0.$
\end{remark}

\begin{remark}
Gutman et al. in \cite{gavbr} conjectured that for any simple graph $G$, $\mathscr{E}(G) \leq \mathscr{LE}(G)$. It was disproved in literature
by showing many counterexamples. However, we confirmed in \cite{dll} that it is almost true by showing that for almost all graphs
$G_p\in \mathcal{G}_{n, p}$, $\mathscr{E}(G_p) < \mathscr{LE}(G_p)$. Next we will give a stronger result.
\end{remark}

For the matrices with degree-distance-based entries, from Theorems \ref{main2} and \ref{dis} we can get the following result.

\begin{theorem}\label{lap-adj}
Let $f(D(i,j), d_i, d_j)$ be real function  symmetric in $i$ and $j$ satisfying that $f(1,(1+o(1))np, (1+o(1))np)=(1+o(1))f(1, np, np)$ and
$f(2,(1+o(1))np, (1+o(1))np)=(1+o(1))f(2, np, np)$. Then for almost all graphs $G_p\in \mathcal{G}_{n, p}$,
$$\mathscr{E}(W_f(G_p)) < \mathscr{LE}_f(G_p),$$
which implies $\mathscr{E}(A_f(G_p)) < \mathscr{LE}_f(G_p)$, and hence, $\mathscr{E}(G_p) < \mathscr{LE}(G_p)$.
\end{theorem}
The result holds true because in Theorems \ref{main2} and \ref{dis} the coefficients $\frac{2\sqrt{2}} {3} > \frac {8} {3\pi}$.

\section[The energies $\mathscr{LEL}_f$ and $\mathscr{IE}_f$]{The energies $\mathscr{LEL}_f$ and $\mathscr{IE}_f$}

\subsection[Weighted adjacency matrices of random graphs]{Weighted adjacency matrices of random graphs}

In this section, we deal with a particular kind of weighted-distance matrix, the weighted adjacency matrix, before handling the general cases. A \emph{weighted adjacency matrix} is a weighted-distance matrix in which $f(D(i, j), d_{i}, d_{j})=0$ if $D(i, j)>1$. We use $A_f(G)$ rather than $W_f(G)$ to represent it, and keep other notations unchanged.

First of all, we would like to introduce some properties of random graphs as preliminaries of our work. As is shown in Section 1, the matrix $W_f(G_p)$ may be rather complicated since the diameter of a graph can be very large, and also many different values of distance $D(i, j)$ for pairs of vertices of $G$ are involved. However, things are not that disappointed from the probability point of view. In fact, there is a famous result, stated as follows.
\begin{lemma}\label{diameter} \cite{boll}
Almost all graphs have diameter two.
\end{lemma}

So, from this result, one can see that to study the asymptotic property, it suffices to deal with graphs of diameter $2$, whose weighted distance matrix $W_f(G_p)$ consists of entries with only values $0$, $f(1, d_{i}, d_{j})$ and $f(2, d_{i}, d_{j})$.

Another fact we need to point out is that for almost all graphs $G_p \in \mathcal{G}_{n, p}$, the vertex-degrees of $G_p$ concentrate around the value $np$.
\begin{lemma} (see \cite{boll}) \label{lem2.2}
Let $\varepsilon>0$ be fixed, and $\varepsilon n^{-3/2}\leq p\leq 1-\varepsilon n^{-3/2}$. Let $q=q(n)$ be a natural number and set\\
\begin{displaymath}
\mu_{q}=nB(q;n-1, p)\quad and\quad\nu_{q}=n\{1-B(q+1;n-1, p)\},
\end{displaymath}
where\\
\begin{displaymath}
B(l;m, p)=\sum\limits_{j\geq l}b(j;m, p)
\end{displaymath}
in which $b(j;m, p)=\binom{m}{j}p^{j}(1-p)^{m-j}$ is subject to the binomial distribution. For a random graph $G\in\mathcal{G}_{n, p}$,
denote by $Y_{q}(G)$ the number of vertices with degrees at least $q$ and $Z_{q}(G)$ the number of vertices with degrees at most $q$. Then\\
\begin{displaymath}
(i) \ \  if\,\,\,\mu_{q}\rightarrow 0,\,P(Y_{q}=0)\rightarrow 0;\, \
(ii) \ \ if\,\,\,\nu_{q}\rightarrow 0,\,P(Z_{q}=0)\rightarrow 0.
\end{displaymath}
\end{lemma}

Remember that $p\in (0,1)$ is a constant. It is not difficult to check that the minimum and maximum degrees $\delta$ and $\Delta$
of a random graph $G_p$ on $n$ vertices satisfy that
\begin{equation}
np-n^{\frac{3}{4}}<\delta(G_p)\leq\Delta(G_p)<np+n^{\frac{3}{4}},\mbox{ a.s.}
\end{equation}
$(i)$ and $(ii)$ hold by Chernoff's Inequality. In summary, we just need to deal with graphs with diameter 2, in which all vertex-degrees fall in the interval $(np-n^{\frac{3}{4}},np+n^{\frac{3}{4}})$.

From now on, we start the main procedure of our proof. Let $L_f^{\prime}=\frac{L_f}{f(np, np)}=\sum\limits_{i\sim j} Q_{ij}$, where $Q_{ij}$ is an $n\times n$ matrix with only four non-zero entries: $Q_{ij}(i, i)=Q_{ij}(j, j)=a_{ij}$, and $Q_{ij}(i, j)=Q_{ij}(j, i)=-a_{ij}$, where every $a_{ij}=1+o(1)>0$. $Q_{ij}$ is semi-positive definite, and so is $L_f^{\prime}$.

Remember that $L(G)=D(G)-A(G)$ is the Laplacian matrix (not weighted). Take $X=(x_1, \ldots x_n)\in \mathbf{R}^n$, $X^{\top} L_f^{\prime}X=X^{\top}(\sum\limits_{i\sim j} Q_{ij})X=\sum\limits_{i\sim j} a_{ij}(x_i-x_j)^2$, and so $rank(L_f^{\prime})=rank(L)=n-1$. Assume that the spectrum of $L$ is $\lambda_1(L)\geq\lambda_2(L)\geq\ldots\geq\lambda_{n-1}(L)>\lambda_n(L)=0$, and the spectrum of $L_f^{\prime}$ is $\lambda_1(L_f^{\prime})\geq\lambda_2(L_f^{\prime})\geq\ldots\lambda_{n-1}(L_f^{\prime})>\lambda_n(L_f^{\prime})=0$. From here we can say that the definitions of $\mathscr{LEL}_f$ and $\mathscr{IE}_f$ are quite reasonable, since the spectrum of $L_f(G)$ is semi-positive or semi-negative definite, depending on the sign of $f(np, np)$.

We regard $L_f^\prime$ as a perturbed matrix of $L$ and will show that the eigenvalues of $L_f^\prime$ approximate those of $L$. The main tool we count on is Weyl's Inequality, which performs well in this kind of problems.
\begin{lemma}\label{weyl}
(Weyl's Inequality, see \cite{weyl}) Assume $M=H+P$, where $M, H$ and $P$ are $n\times n$ Hermitian matrices. Suppose that $M$ has eigenvalues
\begin{displaymath}
\xi_1\geq\xi_2\geq\ldots\geq\xi_n,
\end{displaymath}
and $H$ has eigenvalues
\begin{displaymath}
\nu_1\geq\nu_2\geq\ldots\geq\nu_n,
\end{displaymath}
and $P$ has eigenvalues
\begin{displaymath}
\theta_1\geq\theta_2\geq\ldots\geq\theta_n.
\end{displaymath}
If $j+k-n\geq i\geq r+s-1$, then we have
\begin{displaymath}
\nu_j+\theta_k\leq\xi_i\leq\nu_r+\theta_s.
\end{displaymath}
Especially, $\nu_i+\theta_n\leq\xi_i\leq\nu_i+\theta_1$.
\end{lemma}

Also, we need some information about the eigenvalues of Wigner and Markov matrices to characterize the spectrum of $L$.
\begin{lemma} \label{marad} (Bryc \cite{bryc}). Suppose $M$ is a Markov matrix whose off-diagonal entries are i. i. d. random variables with mean zero, variance one and finite fourth moment. Then the spectrum radius $r(M)\sim\sqrt{2n\log n}$ a.s.
\end{lemma}

\begin{lemma} \label{wigrad} (F$\ddot{u}$redi \emph{et. al.} \cite{fure}) Let $A$ be an $n\times n$ Wigner matrix with $\mathbf{E}A(1, 2)>0$. Then,  with probability $1-o(1)$, the eigenvalues of $A$ are $O(n^{1/2})$ except the largest one, which is $O(n)$.
\end{lemma}

\begin{theorem}
For almost all graphs $G_p\in \mathcal{G}_{n, p}$, $\lambda_i(L_f^{\prime})=(p+o(1))n,$ $i=1, 2, \ldots, n-1$.
\begin{proof} It is easy to see that
$\bar{L}=(p(1-p))^{-1/2}[(D-p(n-1)I)-(A-p(J-I))]=(p(1-p))^{-1/2}(L+pJ-npI)$ is a Markov matrix satisfying the conditions in Lemma \ref{marad}, where $J$ is the $n\times n$ matrix with all ones, and $r(\bar{L})=o(n)$. The eigenvalues of $pJ-npI$ are $-pn$ with multiplicity n-1 and 0 with multiplicity 1. Apply Weyl's Inequality on $\bar{L}$ and $pJ-npI$, we get $\lambda_i(L)=(p+o(1))n$, $i=1, \ldots, n-1$.\\
Decompose $L_f^{\prime}=L+\tilde{L}$, $\tilde{L}=\sum\limits_{i\sim j} \tilde{Q}_{ij}$, $\tilde{Q}_{ij}(i, i)=\tilde{Q}_{ij}(j, j)=\tilde{a}_{ij}$, and $\tilde{Q}_{ij}(i, j)=\tilde{Q}_{ij}(j, i)=-\tilde{a}_{ij}$, where $\tilde{a}_{ij}=a_{ij}-1=o(1)$. Assume that the spectrum of $\tilde{L}$ is $\sigma_1\geq\sigma_2\geq\ldots\geq\sigma_n$. By Rayleigh's Inequality, $\sigma_1=\max\limits_{||X||=1} X^{\top}\tilde{L}X$, $\sigma_n=\min\limits_{||X||=1} X^{\top}\tilde{L}X$. So,  $\sigma_1=o(\lambda_1(L))$ and $\sigma_n=o(\lambda_1(L))$, and thus the spectrum radius of $\tilde{L}$ $r(\tilde{L})=o(\lambda_1(L))$. Apply Weyl's Inequality once again, the proof is thus complete.
\end{proof}
\end{theorem}

Now consider the $f$-weighted signless Laplacian matrix $L_f^+$. Noticing  that $L^+=L+2A$, and using Lemma \ref{wigrad} and Weyl's Inequality, we can get that
\begin{displaymath}
\lambda_1(L^+)=O(n), \lambda_n(L^+)=O(n),\quad and\quad \lambda_i(L^+)=(p+o(1))n\quad for \quad i=2, \ldots, n-1.
\end{displaymath}
Repeating the former process (almost the same and we omit it), we have that
\begin{displaymath}
\lambda_1(L_f^+)=O(nf(np, np)), \lambda_n(L_f^+)=O(nf(np, np)),
\end{displaymath}
and
\begin{displaymath}
\lambda_i(L_f^+)=f(np, np)(p+o(1))n\quad for\quad i=2, \ldots, n-1.
\end{displaymath}
Then we can obtain the values of $\mathscr{LEL}_f$ and $\mathscr{IE}_f$
as follows.
\begin{corollary} \label{vice}
Let $f(d_i, d_j)$ be a real function symmetric on $i$ and $j$ with the property that $f((1+o(1))np, (1+o(1))np)=(1+o(1))f(np, np)$.
Then, for almost all graphs $G_p\in \mathcal{G}_{n, p}$,
\begin{displaymath}
\mathscr{LEL}_f(G_p)=\sqrt{|f(np, np)|}(\sqrt{p}+o(1))\cdot n^{3/2},
\end{displaymath}
\begin{displaymath}
\mathscr{IE}_f(G_p)=\sqrt{|f(np, np)|}(\sqrt{p}+o(1))\cdot n^{3/2}.
\end{displaymath}
\end{corollary}

\subsection[General case]{General case}

In this subsection, we deal with the general case. In chemical practice, $f(1, np, np)$, $f(2, np, np)$ are often functions with a definite limit as $n\rightarrow\infty$. In fact, the indices already defined so far are all positive power functions. So, we suppose that $f(1, np, np)/f(2, np, np)\rightarrow C$ (may be $\infty$).

First of all, consider the special case $f(1, np, np)\sim f(2, np, np)$. Then $\frac{L_f}{f(1, np, np)}$ can be decomposed into the sum of two matrices,
one is the Laplacian matrix of the complete graph $K_n$ and the other is the perturbed matrix whose off-diagonal entries are all $o(1)$.
Using the method in the last section and noticing the fact that the Laplacian eigenvalues of $K_n$ are $n$ with multiplicity $n-1$ and 0 with multiplicity 1,
and the signless Laplacian eigenvalues of $K_n$ are $n-2$ with multiplicity $n-1$ and $2n-2$ with multiplicity 1, we get that $\mathscr{LEL}_f=\sqrt{|f(2, np, np)|}(1+o(1))n^{3/2}$
and $\mathscr{IE}_f=\sqrt{|f(2, np, np)|}(1+o(1))n^{3/2}$.

From now on, we reckon $f(1, np, np)/f(2, np, np)\nrightarrow 1$. Decompose $L_f$ into two matrices $L_f=L_1+L_2$, where
\begin{displaymath}
L_{1}(i, j)=
\begin{cases}
-(f(1, d_{i}, d_{j})-f(2, d_{i}, d_{j})),&\quad \mbox{$i$ and $j$ are adjacent,}\\
0,&\quad \mbox{$i$ and $j$ are nonadjacent,}\\
-\sum\limits_{k\neq i} L_1(i, k),&\quad \mbox{$i=j$},
\end{cases}
\end{displaymath}
and
\begin{displaymath}
L_2(i, j)=
\begin{cases}
-f(2, d_{i}, d_{j}),&\quad \mbox{$i\neq j$,}\\
-\sum\limits_{k\neq i} L_2(i, k),&\quad \mbox{$i=j$.}
\end{cases}
\end{displaymath}

For convenience, denote $F(d_i, d_j)=f(1, d_i, d_j)-f(2, d_i, d_j)$, with the property that $F((1+o(1))np, (1+o(1))np)=(1+o(1))F(np, np)$.
According to the proof in the last section, the spectrum of $\frac{L_1}{F(np, np)}$ is $u_1^\prime\geq\ldots\geq u_{n-1}^\prime >u_n^\prime=0$,
with $u_i^\prime=(p+o(1))n$ for $1\leq i<n$; and the spectrum of $\frac{L_2}{f(2, np, np)}$ is $v_1^\prime\geq\ldots\geq v_{n-1}^\prime >v_n^\prime=0$,
with $v_i^\prime=(1+o(1))n$ for $1\leq i<n$. As the signs of $F(np, np)$ and $f(2, np, np)$ determine the spectrum of $L_1$ and $L_2$,
we distinguish four cases to discuss, in each of which the Weyl's Inequality is applied.

{\bf Case 1.} $F(np, np)>0$, $f(2, np, np)\geq 0$.

Then, $\lambda_i(L_f)=F(np, np)(p+o(1))n+f(2, np, np)(1+o(1))n$ for $i=1, \ldots, n-2$, $\lambda_{n-1}=O(nF(np, np))$.

{\bf Case 2.} $F(np, np)>0$, $f(2, np, np)<0$.

Then, $\lambda_1(L_f)=O(nF(np, np))$, $\lambda_n(L_f)=O(nF(np, np))$ and $\lambda_i(L_f)=F(np, np)(p+o(1))n+f(2, np, np)(1+o(1))n$ for $i=2, \ldots, n-1$.

{\bf Case 3.} $F(np, np)<0$, $f(2, np, np)\geq 0$.

Then, $\lambda_{1}(L_f)=O(nF(np, np))$, $\lambda_n(L_f)=O(nF(np, np))$ and $\lambda_i(L_f)=F(np, np)(p+o(1))n+f(2, np, np)(1+o(1))n$ for $i=2, \ldots, n-1$.

{\bf Case 4.} $F(np, np)<0$, $f(2, np, np)< 0$.

Then, $\lambda_{n}(L_f)=O(nF(np, np))$ and $\lambda_i(L_f)=F(np, np)(p+o(1))n+f(2, np, np)(1+o(1))n$ for $i=2, \ldots, n-1$.

For the signless Laplacian matrix $L_f^+$, the method is the same. Let $L_f^+=L_1^++L_2^+$, where
\begin{displaymath}
L_1^+(i, j)=
\begin{cases}
f(1, d_{i}, d_{j})-f(2, d_{i}, d_{j}),&\quad \mbox{$i$ and $j$ are adjacent,}\\
0,&\quad \mbox{$i$ and $j$ are nonadjacent,}\\
\sum\limits_{k\neq i} L_1^+(i, k),&\quad \mbox{$i=j$},
\end{cases}
\end{displaymath}
and
\begin{displaymath}
L_2^+(i, j)=
\begin{cases}
f(2, d_{i}, d_{j}),&\quad \mbox{$i\neq j$,}\\
\sum\limits_{k\neq i} L_2^+(i, k),&\quad \mbox{$i=j$.}
\end{cases}
\end{displaymath}
The spectrum of $\frac{L_1}{F(np, np)}$ is $u_1^{\prime+}\geq\ldots\geq u_{n-1}^{\prime+} \geq u_n^{\prime+}$,
with $u_1^{\prime+}=O(n)$, $u_n^{\prime+}=O(n)$, and $u_i^{\prime+}=(p+o(1))n$ for $1<i<n$.
The spectrum of $\frac{L_2}{f(2, np, np)}$ is $v_1^{\prime+}\geq\ldots\geq v_{n-1}^{\prime+} \geq v_n^{\prime+}$,
with $v_1^{\prime+}=(2+o(1))n$ and $v_i^{\prime+}=(1+o(1))n$ for $1<i\leq n$. The following two cases are distinguished.

{\bf Case 1.} $f(2, np, np)\geq 0$.

Then, $\lambda_1(L_f^+)=O(nF(np, np))$, $\lambda_2(L_f^+)=O(nF(np, np))$,
$\lambda_n(L_f^+)=O(nF(np, np))$, and $\lambda_i(L_f^+)=F(np, np)(p+o(1))n+f(2, np, np)(1+o(1))n$ for $i=3, \ldots, n-1$.

{\bf Case 2.} $f(2, np, np)<0$.

Then, $\lambda_1(L_f^+)=O(nF(np, np))$, $\lambda_{n-1}(L_f^+)=O(nF(np, np))$, $\lambda_n(L_f^+)=O(nF(np, np))$
and $\lambda_i(L_f^+)=F(np, np)(p+o(1))n+f(2, np, np)(1+o(1))n$ for $i=2, \ldots, n-2$.

From the discussion above, we can see that almost all eigenvalues of $L_f$ ($L_f^+$) are equal to
$F(np, np)(p+o(1))n+f(2, np, np)(1+o(1))n$. So, we can calculate $\mathscr{LEL}_f$ and $\mathscr{IE}_f$
directly by discussing two cases: $C=\pm\infty$ and $C\in \mathbf{R}$. Finally, we can get Theorem \ref{main1}.

\section[(Signless) Laplacian energy of distance-based matrices]{(Signless) Laplacian energy of distance-based matrices}

On the whole, the calculation of (signless) Laplacian energy of weighted distance matrices, or even degree-weighted adjacency
matrices are much more complicated than $\mathscr{LEL}_f$ and $\mathscr{IE}_f$. We can only give a magnitude at present,
which is far from satisfactory. However, if we restrict the study to distance-based matrices, that is, $f(D(i, j), d_i, d_j)$
is a function depending only on $D(i, j)$, nice lower and upper bounds can be neatly obtained.

Denote $f(1, d_i, d_j)=D_1$ and $f(2, d_i, d_j)=D_2$. We continue to adopt the notations in Section 2,
but replace $L_f$ with $L_D$ to emphasize the object we are studying. Let $L_D=L_1+L_2$, here $L_1=(D_1-D_2)L$
and $L_2=D_2L_0$, where $L_0$ is the Laplacian matrix of $K_n$. The eigenvalues of $L_2$ are $D_2n$ with multiplicity $n-1$ and 0 with multiplicity 1.
It is clear that we only need to consider the case that $D_1-D_2> 0$; otherwise, replace $L_D$ with $-L_D$. So, we set it as a premise of this section.

First, assume $D_2\geq 0$. Using Weyl's Inequality, we get
\begin{displaymath}
(D_1-D_2)\lambda_{i+1}(L)+D_2n\leq \lambda_i(L_D)\leq  (D_1-D_2)\lambda_{i}(L)+D_2n,\quad i=1,\ldots n-1,
\end{displaymath}
and $\lambda_n(L_D)=0$, from which we have
\begin{displaymath}
(D_1-D_2)\sum\limits_{i=2}^n\lambda_i(L)+D_2n(n-1)\leq\sum\limits_{i=1}^{n-1}\lambda_i(L_D)\leq (D_1-D_2)\sum\limits_{i=1}^{n-1}\lambda_i(L)+D_2n(n-1).
\end{displaymath}

We use $\bar{\lambda}$ to represent the average value of the eigenvalues. That is,
\begin{displaymath}
\bar{\lambda}(L_D)=(D_1-D_2)\bar{\lambda}(L)+D_2n+O(1).
\end{displaymath}
Let $\Lambda=(D_1-D_2)\bar{\lambda}(L)+D_2n$, and let $s$ be the integer such that $\lambda_s(L)\geq \bar{\lambda}(L)>\lambda_{s+1}(L)$. Then,
\begin{align}
\mathscr{LE}_D(G_p)&=\sum\limits_{i=1}^{n} \mid\lambda_i(L_D)-\frac{\sum\limits_{i\neq j}f(D(i, j), d_i, d_j)}{n}\mid=\sum\limits_{i=1}^{n} \mid\lambda_i(L_D)-\bar{\lambda}(L_D)\mid\\
&=\sum\limits_{i=1}^{n} \mid\lambda_i(L_D)-\Lambda\mid+O(n)=\sum\limits_{i=1}^{s-1} \mid\lambda_i(L_D)-\Lambda\mid+\sum\limits_{i=s+1}^{n-1} \mid\lambda_i(L_D)-\Lambda\mid+O(n).
\end{align}
Noticing that $(D_1-D_2)(\lambda_{i+1}(L)-\bar{\lambda}(L))\leq \lambda_i(L_D)-\Lambda\leq (D_1-D_2)(\lambda_{i}(L)-\bar{\lambda}(L))$ for $i=1,\ldots n-1$, we get
\begin{displaymath}
\sum\limits_{i=2}^s (D_1-D_2)(\lambda_{i}(L)-\bar{\lambda}(L))\leq\sum\limits_{i=1}^{s-1} \mid\lambda_i(L_D)-\Lambda\mid\leq\sum\limits_{i=1}^{s-1} (D_1-D_2)(\lambda_{i}(L)-\bar{\lambda}(L)),
\end{displaymath}
\begin{displaymath}
\sum\limits_{i=s+1}^{n-1} (D_1-D_2)(\bar{\lambda}(L)-\lambda_i(L))\leq\sum\limits_{i=s+1}^{n-1} \mid\lambda_i(L_D)-\Lambda\mid\leq\sum\limits_{i=s+2}^{n} (D_1-D_2)(\bar{\lambda}(L)-\lambda_i(L)).
\end{displaymath}
Thus,
\begin{displaymath}
\sum\limits_{i=1}^n \mid\lambda_i(L_D)-\Lambda\mid=(D_1-D_2)\sum\limits_{i=1}^n \mid\lambda_i(L)-\bar{\lambda}(L)\mid+O(n).
\end{displaymath}
Recall from Theorem \ref{lap} that the magnitude of the asymptotic Laplacian energy is $O(n^{\frac{3}{2}})$. Thus, $\mathscr{LE}_D(G_p)$ is essentially $|D_1-D_2|$ times $\mathscr{LE}(G_p)$. As for the case $D_2<0$, the interlacing relationship between $\lambda_i(L_D)$ and $\lambda_i(L)$ becomes
\begin{displaymath}
(D_1-D_2)\lambda_i(L)+D_2n\leq\lambda_i(L_D)\leq (D_1-D_2)\lambda_{i-1}(L)+D_2n, \quad i=2,\ldots,n.
\end{displaymath}
Repeating the process above, the identical conclusion can be reached.

$\mathscr{LE}_f^+(G_p)$ can be calculated in the same way. Let $L_f^+=L_1^++L_2^+$, where $L_1^+=(D_1-D_2)L^+$ and $L_2^+=D_2L_0^+$. The eigenvalues of $L_2^+$ are $2D_2(n-1)$ with multiplicity 1 and $D_2(n-2)$ with multiplicity $n-1$. Using Weyl's Inequality, the two relations can be obtained, respectively, as follows.

If $D_2\geq 0$, then $\lambda_1^+(L_D)=O(n)$ and
\begin{displaymath}
(D_1-D_2)\lambda_i^+(L)+D_2(n-2)\leq\lambda_i^+(L_D)\leq(D_1-D_2)\lambda_{i-1}^+(L)+D_2(n-2),\quad i=2,\ldots,n.
\end{displaymath}

If $D_2<0$, then $\lambda_n^+(L_D)=O(n)$ and
\begin{displaymath}
(D_1-D_2)\lambda_{i+1}^+(L)+D_2(n-2)\leq\lambda_i^+(L_D)\leq(D_1-D_2)\lambda_{i}^+(L)+D_2(n-2),\quad i=1,\ldots,n-1.
\end{displaymath}
Similarly, $\mathscr{LE}_D^+$ is essentially $|D_1-D_2|\mathscr{LE}_D^+$.
Combining Theorems \ref{lap} and \ref{sign}, we eventually get the asymptotic value of the (signless) Laplacian energy, which is stated in Theorem \ref{main2}.

\section[Applications for matrices with distance-based and degree-distance-based weights of chemical use]{Applications for matrices with distance-based and degree-distance-based weights of chemical use}

As once proposed in \cite{XL} and commented in \cite{llw}, if we use a matrix to represent the structure of a molecular graph with weights separately on its pairs of vertices, it will completely keep the structural information of the graph, i.e., a matrix keeps much more structural information than a numeral of an index. So, it is significant to study the properties of these chemical-index-weighted matrices. Here we use three tables to list our results, in terms of three different kinds of indices, that is, the degree-based only \cite{lls}, distance-based only, and the degree-distance-based mixed \cite{llw}. One can easily calculate these results from the formulas given in the above sections, details of which are omitted.

\begin{table}[H]
\caption{Degree-based indices}
\begin{tabular}{ccc}
\toprule
  Index & $f(d_i, d_j)$ & $\mathscr{LEL}_f$ and $\mathscr{IE}_f$\\
  \midrule
  First Zagreb & $d_i+d_j$ & $(\sqrt{2}p+o(1))\cdot n^2$\\
  Second Zagreb & $d_id_j$ & $(p\sqrt{p}+o(1))\cdot n^{5/2}$\\
  Randi$\acute{c}$ & $\frac{1}{\sqrt{d_id_j}}$ & $(1+o(1))n$\\
  General Randi$\acute{c}$ & $(d_id_j)^\alpha$ & $(p^{\alpha+\frac{1}{2}}+o(1))\cdot n^{\alpha+\frac{3}{2}}$\\
  ABC & $\frac{\sqrt{d_i+d_j-2}}{\sqrt{d_id_j}}$ & $((2p)^{\frac{1}{4}}+o(1))\cdot n^{5/4}$\\
  AZI & $(\frac{d_id_j}{d_i+d_j-2})^3$ & $(\frac{p^2}{2\sqrt{2}}+o(1))\cdot n^3$\\
  AG & $\frac{2\sqrt{d_id_j}}{d_i+d_j}$ & $(\sqrt{p}+o(1))\cdot n^{3/2}$\\
  Harmonic & $\frac{2}{d_i+d_j}$ & $(1+o(1))n$\\
  SCI & $\frac{1}{\sqrt{d_i+d_j}}$ & $((\frac{p}{2})^{\frac{1}{4}}+o(1))\cdot n^{5/4}$\\
  First multiple Zagreb & $\frac{\log d_i}{d_i}+\frac{\log d_j}{d_j}$ & $(\sqrt{2}+o(1))n\sqrt{\log n}$\\
  Modified multiple Zagreb & $\log(d_i+d_j)$ & $(\sqrt{p}+o(1))\cdot n^{3/2}\sqrt{\log n}$\\
  Second multiple Zagreb & $\log d_i+\log d_j$ & $(\sqrt{2p}+o(1))\cdot n^{3/2}\sqrt{\log n}$\\
  Lanzhou & $(n-1)(d_i+d_j)-(d_i^2+d_j^2)$ & $(p\sqrt{2(1-p)}+o(1))\cdot n^{5/2}$\\
  \bottomrule
  \end{tabular}
  \end{table}

  \begin{table}[H]
  \caption{Distance-based indices}
  \begin{tabular}{ccc}
  \toprule
  Index & $f(D(i, j))$ & $\mathscr{LEL}_f$ and $\mathscr{IE}_f$\\
  \midrule
  Harary \cite{plav, ivanc} & $1/D(i, j)$ & $(\sqrt{\frac{1+p}{2}}+o(1))\cdot n^{3/2}$\\
  Hyper-Wiener & $\frac{1}{2}(D(i, j)+D^{2}(i, j))$ & $(\sqrt{3-2p}+o(1))\cdot n^{3/2}$\\
  Reciprocal complementary Wiener \cite{iib} & $\frac{1}{diam(G)^{\rm *}+1-D(i, j)}$ & $(\sqrt{1-\frac{p}{2}}+o(1))\cdot n^{3/2}$\\
  Reverse Wiener \cite{van} & $diam(G)-D(i, j)$ & $(\sqrt{p}+o(1))\cdot n^{3/2}$\\
  \bottomrule
  \end{tabular}
  \footnotesize{$^{\rm *}diam(G)$ is the diameter of graph $G$.}
  \end{table}

  \begin{table}[H]
  \caption{Degree-distance-based indices}
  \begin{tabular}{ccc}
  \toprule
  Index & $f(D(i, j), d_i, d_j)$ & $\mathscr{LEL}_f$ and $\mathscr{IE}_f$\\
  \midrule
  degree-distance-index \cite{dob} & $(d_{i}+d_{j})D(i, j)$ & $(\sqrt{4p-2p^2}+o(1))\cdot n^2$\\
  Gutman \cite{gut2} & $d_{i}d_{j}D(i, j)$ & $(\sqrt{2p^2-p^3}+o(1))\cdot n^{5/2}$\\
  Additively weighted Harary \cite{ali, hua} & $\frac{d(u)+d(v)}{D(u, v)}$ & $(\sqrt{p+p^2}+o(1))\cdot n^2$\\
  Multiplicatively weighted Harary \cite{ali, hua} & $\frac{d_{i}d_{j}}{D(i, j)}$ & $(\sqrt{\frac{p^2+p^3}{2}}+o(1))\cdot n^{5/2}$\\
  \bottomrule
  \end{tabular}
  \end{table}

\begin{remark} From the above three tables, one can see the following facts. (1) \ It seems that some values, say $(1+o(1))n$, do not depend on the probability $p$.
However within the term $o(1)$, $p$ will plays a role, though not in the first domination term. (2) \ Most of the values are monotonic when $p$ gets large.
But, some of the values get their extremal at nontrivial points, say the Gutman index.
\end{remark}

\end{document}